\documentclass{article}
\usepackage{graphicx} % Required for inserting images
\usepackage{amsmath,amssymb,amsfonts,amsthm}
\usepackage{soul}
\usepackage{mathtools}
\usepackage{dsfont}
\usepackage{booktabs}
\usepackage{caption}
\usepackage{mathabx}
\usepackage{tablefootnote}
\usepackage{titlesec}
\usepackage{algorithmicx}
\usepackage{bm}
\usepackage{setspace}
\usepackage{hyperref}
\usepackage{algorithm}
\usepackage{array}
\usepackage{placeins}
\usepackage{algpseudocode}
\usepackage{xspace}
\usepackage{textgreek}
\usepackage{enumitem}
\makeatletter
\renewcommand{\ALG@beginalgorithmic}{\footnotesize}
\makeatother
\usepackage{xcolor}
\usepackage{natbib}
\usepackage{tikz}
\usepackage{url}
\usepackage{threeparttable}
\usepackage{listings}
\usetikzlibrary{calc,positioning,arrows.meta,quotes,decorations.pathreplacing}
\usepackage{pgfplots}
\pgfplotsset{compat=1.17}
\newcommand{\ourmethod}{HALLaR\xspace}
\newcommand{\gpuourmethod}{cuHALLaR\xspace}

\usepackage{amsmath}
\usepackage{amsfonts}
\usepackage{latexsym}
\usepackage{amssymb}
\usepackage{stmaryrd}
\usepackage{mathtools}

\newtheorem{thm}{Theorem}[section]

\newtheorem{remark}[thm]{Remark}

\newcommand{\beq}{\begin{equation}}
\newcommand{\eeq}{\end{equation}}
\newcommand{\beqa}{\begin{eqnarray}}
\newcommand{\eeqa}{\end{eqnarray}}
\newcommand{\beqas}{\begin{eqnarray*}}
\newcommand{\eeqas}{\end{eqnarray*}}
\newcommand{\ei}{\end{itemize}}

\def\defi{\vcentcolon=}

\newcommand{\R}{\mathbb{R}}

\newcommand{\ignore}[1]{}

\newtheorem{innercustomgeneric}{\customgenericname}
\providecommand{\customgenericname}{}
\newcommand{\newcustomtheorem}[2]{%
  \newenvironment{#1}[1]
  {%
   \renewcommand\customgenericname{#2}%
   \renewcommand\theinnercustomgeneric{##1}%
   \innercustomgeneric
  }
  {\endinnercustomgeneric}
}
\newcustomtheorem{blackbox}{Blackbox}

\allowdisplaybreaks

\title{An Accelerated Hybrid Low-Rank Augmented Lagrangian Method for Large-Scale Semidefinite Programming on GPU}

\author{
{Jacob M. Aguirre,}
\thanks{H. M. Stewart School of Industrial and Systems Engineering, Georgia Tech, Atlanta, GA, 30332-0205. (Email: {\tt aguirre@gatech.edu})}
{\quad Diego Cifuentes,}
\thanks{H. M. Stewart School of Industrial and Systems Engineering, Georgia Tech, Atlanta, GA, 30332-0205. (Email: {\tt dfc3@gatech.edu})}
{\quad Vincent Guigues}
\thanks{School of Applied Mathematics, FGV,
Praia de Botafogo, Rio de Janeiro, Brazil. (Email:{\tt vincent.guigues@fgv.br})}
{Renato D.C. Monteiro,}
\thanks{H. M. Stewart School of Industrial and Systems Engineering, Georgia Tech, Atlanta, GA, 30332-0205. (Email: {\tt rm88@gatech.edu})}
{\quad Victor Hugo Nascimento,}
\thanks{School of Applied Mathematics, FGV,
Praia de Botafogo, Rio de Janeiro, Brazil. (Email:{\tt nascimento.victor.1@fgv.edu.br})}
{\quad Arnesh Sujanani.}
\thanks{Department of Combinatorics and Optimization, University of Waterloo, Waterloo, ON, N2L 3G1.(Email: {\tt a3sujana@uwaterloo.edu}) \\~\\~\\~\\ Authors are listed in alphabetical order}
}

\date{\today}

\title{
A User Manual for \gpuourmethod: A GPU Accelerated Low-Rank Semidefinite Programming Solver
\thanks{\textbf{Funding}: 
Jacob M. Aguirre is supported by the National Science Foundation Graduate Research Fellowship under Grant No. DGE-2039655.
Diego Cifuentes is supported by U.S. Office of Naval Research, N00014-23-1-2631.
Renato D.C. Monteiro is supported by AFOSR Grant FA9550-25-1-0131. 
}}

\begin{document}
\author{Jacob M. Aguirre\thanks{H. M. Stewart School of Industrial and Systems Engineering, Georgia Tech, Atlanta, GA, 30332-0205. \protect\protect\href{mailto:aguirre@gatech.edu}{aguirre@gatech.edu}}\hspace*{0.5em}
\and Diego Cifuentes\thanks{H. M. Stewart School of Industrial and Systems Engineering, Georgia Tech, Atlanta, GA, 30332-0205.  
\protect\protect\href{mailto:dfc3@gatech.edu}{dfc3@gatech.edu}}
\and Vincent Guigues\thanks{School of Applied Mathematics, FGV,
Praia de Botafogo, Rio de Janeiro, Brazil. 
\protect\protect\href{mailto:vincent.guigues@fgv.br}{vincent.guigues@fgv.br}}
\and Renato D.C. Monteiro\thanks{School of Industrial and Systems Engineering, Georgia Institute of
Technology, Atlanta, GA, 30332-0205. 
\protect\protect\href{mailto:rm88@gatech.edu}{rm88@gatech.edu}}
\and Victor Hugo Nascimento\thanks{School of Applied Mathematics, FGV,
Praia de Botafogo, Rio de Janeiro, Brazil.
\protect\protect\href{mailto:vhn@fgv.br}{vhn@fgv.br}}
\and Arnesh Sujanani\thanks{Department of Combinatorics and Optimization, University of Waterloo, Waterloo, ON, N2L 3G1.
\protect\protect\href{mailto:a3sujana@uwaterloo.ca}{a3sujana@uwaterloo.ca}\\~\\ Authors are listed in alphabetical order.}
}

\setcitestyle{square}
\maketitle

\begin{abstract}
We present a Julia-based interface to the precompiled \ourmethod and \gpuourmethod binaries for large-scale semidefinite programs (SDPs). Both solvers are established as fast and numerically stable, and accept problem data in formats compatible with SDPA and a new enhanced data format taking advantage of Hybrid Sparse Low-Rank (HSLR) structure. The interface allows users to load custom data files, configure solver options, and execute experiments directly from Julia. A collection of example problems is included, including the SDP relaxations of the Matrix Completion and Maximum Stable Set problems.

\end{abstract}

\section{Introduction}\label{Introduction}
This document serves as a user guide for \ourmethod  \citep{monteiro2024low} and \gpuourmethod  \citep{aguirre2025cuhallar}. Their binaries can be downloaded from \url{https://github.com/OPTHALLaR}. 

Let $\mathbb S^n$ be the set of
$n\times n$ symmetric matrices.
The notation $A \succeq B$ means
that $A-B$ is positive semidefinite.
\gpuourmethod and \ourmethod solve the primal-dual pair of semidefinite programs (SDPs)
\begin{gather*}
\label{eq:sdp-primal}\tag{P}
P_{*}: = \quad \min_{X} \quad \{C \bullet X
\quad : \quad
\mathcal A(X)=b ,\quad
\mathrm{Tr}(X)\leq\tau,\quad X\succeq0 \}
\end{gather*}
and
\begin{gather*}
\label{eq:sdp-dual}\tag{D}
D_{*}: = \quad \max_{(p,S,\theta)}   \quad \{ -b^{\top}p-\tau\theta
\quad : \quad
C+\mathcal A^{*}(p)+\theta I -S =0, \quad S \succeq 0, \quad
\theta \geq 0 \}
\end{gather*}
where $b \in \mathbb R^m$, $C\in \mathbb S^n$, and $\mathcal A: \mathbb S^n \to \mathbb{R}^m$ and $\mathcal A^*: \mathbb{R}^m \to \mathbb S^n$ are linear maps such that 
\begin{align}\label{Definition of Linear Operators}
     \mathcal A(X)= \begin{bmatrix}
           A_1\bullet X \\
           A_2\bullet X \\
           \vdots \\
           A_m \bullet X
         \end{bmatrix}, \quad \mathcal A^{*}(p)=\sum_{i=1}^{m}p_iA_i
        \end{align}
where $A_{\ell} \in \mathbb S^n$ for $\ell=1,\ldots m$. We refer to $\tau$ as the \textbf{trace bound}. We also define
$\Delta^n$ to be the spectraplex, i.e., 
\begin{equation}\label{Delta Definition}
\Delta^n \defi \{ X \in \mathbb S^n : \mathrm{Tr}(X) \leq \tau, X \succeq 0 \}.
\end{equation}

% where $b \in \mathbb{R}^m$, $C\in \mathbb S^n$, $\mathcal A: \mathbb S^n \to \mathbb{R}^m$ is a linear map, and
% $\mathcal A^*: \mathbb{R}^m \to \mathbb S^n$ is its adjoint.
% and $\Delta^n$ is the spectraplex
% \begin{equation}\label{Delta Definition}
% \Delta^n \defi \{ X \in \mathbb S^n : \mathrm{Tr}(X) \leq 1, X \succeq 0 \}.
% \end{equation}

Both solvers are based on an augmented Lagrangian framework with hybrid low‐rank updates, incorporating the Frank--Wolfe method and an adaptive accelerated inexact proximal-point method (ADAP-AIPP) which uses ideas from \cite{Aaronetal2017, WJRproxmet1, WJRComputQPAIPP, CatalystNC, SujananiMonteiro}. \gpuourmethod is a GPU‐accelerated variant intended for high‐throughput computation on modern CUDA‐capable devices, while \ourmethod is a CPU‐based implementation suitable for systems without GPU support. The solvers are numerically stable, memory‐efficient, and support SDPA and Hybrid Sparse Low Rank (HSLR) input formats (discussed more in Section~\ref{sec:interface}).

The guide that follows explains how to prepare problem data in the required HSLR format, how to configure solver parameters either via command‐line options or configuration files, and how to interpret the output files produced.

\section{Installation \& Running}\label{sec:install-run}

\paragraph{Requirements.} \gpuourmethod is distributed as a precompiled binary. No Julia installation is required. A CUDA–enabled NVIDIA GPU, the matching NVIDIA driver, and a compatible CUDA runtime must be available on the system. All Julia dependencies (\texttt{CUDA.jl}, \texttt{LinearAlgebra}, \texttt{SparseArrays}, \texttt{Parameters.jl}, \texttt{KrylovKit.jl}) are bundled. \ourmethod is the CPU version of our code and requires no CUDA libraries to run. We tested \gpuourmethod with CUDA~12.9 on GPUs with compute capability $\ge 5.0$.

\paragraph{Obtain the binary.} Download the precompiled executable from \url{https://github.com/OPTHALLaR}, place it in a working directory, and ensure that the NVIDIA driver and CUDA runtime libraries are discoverable by the dynamic loader.

\paragraph{Quick check.} Verify the installation by running the built–in tests:
\begin{verbatim}
$ ./cuHallar --run_tests
\end{verbatim}
This executes several small example instances shipped with the binary and reports if everything was installed correctly and the code successfully terminated.

\paragraph{Directory layout.} It is not necessary to keep \ourmethod and \gpuourmethod in separate folders. Each folder contains its own executable, default configuration file, and output structure. Invoke the solver from its folder so that the relative paths for the input, output, and configuration files are correctly resolved.

\paragraph{Required inputs.} 

\begin{itemize}
    \item \textbf{config and output files}
    
The config file is a text file that contains the parameter options of \gpuourmethod. The options are organized in the \texttt{key = value} format. For examples of parameter options that can be specified by the user in the text file, the user should refer to Table~\ref{tab:options} in Subsection~\ref{subsec:params}.

% with settings organized , where \texttt{key} is one of the settings from Table~\ref{tab:options}.

The output file is a text file, similar to a CSV file, where a dual solution $p_{*} \in \mathbb R^{m}$ and a low-rank factor $Y_{*} \in \mathbb R^{n\times r}$ of the primal solution $X_{*}$ are saved. We make no assumptions about the naming or the extensions of either the configurations or output file. If extensions such as \texttt{-c}, \texttt{-p}, or \texttt{-d} are omitted, compiled defaults are used. However, the correct files must be provided immediately after each option flag. For instance, if the user passes the wrong file type to the \texttt{-c} option flag, it will get the error listed below.
% the first line of the model file in this example.
\begin{verbatim}
ERROR: ArgParse.ArgParseError("unrecognized option --3 4")
Stacktrace:
 [1] argparse_error(x::Any)
   @ ArgParse ~/.julia/packages/ArgParse/mpp98/src/parsing.jl:9
 [2] parse_long_opt!(state::ArgParse.ParserState, settings::ArgParse.ArgParseSettings)
   @ ArgParse ~/.julia/packages/ArgParse/mpp98/src/parsing.jl:842
 [3] parse_args_unhandled(args_list::Vector, settings::ArgParse.ArgParseSettings, 
    truncated_shopts::Bool)
\end{verbatim}

Here, \texttt{3 4} is not a parameter option for the config file.

\item \textbf{model files: HSLR or SDPA format}

\begin{itemize}
    \item \textbf{HSLR (Hybrid Sparse Low-Rank) format}

    When solving large-scale SDP problem instances with dense cost and constraint matrices, the model size often becomes an issue in terms of storage and RAM. For instance, describing such large SDP instances in standard SDPA format often requires several GBs of storage/RAM. As a result, we introduce the new HSLR format which allows the users to input cost and constraint matrices, $C$ and $A_{\ell}$, $\ell=1,\ldots m$, that are sums of sparse and low-rank matrices. For the sparse component of a matrix, the user only needs to specify the values and indices of the nonzero entries of the upper triangular part of that component. For the low-rank component, the user only needs to input the factors that make-up the low-rank factorization. The sparse and low-rank components are stored internally as CSC (sparse) and dense column–major (low–rank) arrays, respectively.  \\

    HSLR format also requires the user to input dimension pair $(m,n)$, the right-hand side vector $b$, and a trace bound $\tau$ for the trace constraint in \eqref{eq:sdp-primal}. For more details and concrete examples that display how to input these quantities and the cost and constraint matrices correctly for HSLR format, see Subsection~\ref{subsec:data-input}.

    \item \textbf{Sparse SDPA format} 
    
    Standard sparse SDPA \texttt{.dat-s}  files developed in \citep{fujisawa2000sdpa} are also accepted by our \gpuourmethod. 
    % Currently, \gpuourmethod views all SDPA blocks as a single $n\times n$ block diagonal matrix. 
    Each constraint matrix $A_\ell$ that the user provides for SDPA format is imported as a sparse matrix so its low–rank factor is assumed to be zero.
    % and converted in memory to the hybrid layout.
    
As the standard SDPA format lacks support for inputting a trace bound $\tau$ for the trace constraint in \eqref{eq:sdp-primal}, users who wish to use this format \textbf{must} specify the field \texttt{--trace\_bound} in either the command line or in the config file.
% or in the config file, the right-hand side $\tau$ for the trace constraint in \eqref{eq:sdp-primal} by using the command \.
%
% , either in the command-line or in the configuration file. 
%
\gpuourmethod will not run if a user does not provide $\tau$ as an input and the interface will also provide a warning.

% Using \textbf{SDPA format}: \ourmethod also allows the use of the sparse SDPA format. However, the SDPA format does not provide a mechanism to pass the trace bound $\tau$, so users are required to provide one as input.

\end{itemize}

\end{itemize}

\paragraph{Optional input} 

\begin{itemize}
    \item \textbf{initial point}

    An initial primal iterate $Y_0 \in \mathbb R^{n\times r}$ can be supplied by the user as a dense CSV file with no header:
\begin{verbatim}
$ ./cuHallar -i model.hslr -w <path_to_Y0_file>
\end{verbatim}
The initial primal iterate $Y_0$ that the user provides should satisfy $\|Y_0\|^2_{F}\leq \tau$ since this ensures that $\mathrm{Tr}(Y_0Y_0^\top)\le\tau$. If the user does not specify an initial point $Y_0$, then the code internally generates its own $Y_0$ which satisfies $\|Y_0\|_{F}^2\leq \tau$. The user does not need to provide an initial dual iterate $p_0$ since the code always internally sets $p_0$ to be the vector of all zeros.

\end{itemize}

\paragraph{Overriding parameters inline.} Any parameter option of \gpuourmethod may be also supplied on the command line instead of the config file. For example, the following line may be used in the command line
\begin{verbatim}
$ ./cuHallar -i model.hslr --L_inc_fista 3.1 --eps_gap 1e-6 -o out.csv
\end{verbatim}
If a parameter value was specified both on the command line and the config file, the value supplied on the command line takes precedence over the value specified on the config file.

\paragraph{Basic run.}

\begin{itemize}
    \item 
    % input
    % factors of matrices with low-rank  in the form \(A = PDP^\top\), where \(P\) is a \(n \times r\) and \(D\) is a symmetric \(r \times r\), with \(r \ll n\). In this format, the objective and constraints are described as the sum of a sparse matrix and a low-rank factorization. 
    If the user has prepared data in HSLR format and created the required config and output files, the user should then execute the command
\[
\texttt{./cuHallar -i <path\_to\_HSLR\_file> [-c <path\_to\_config\_file>] [-o <path\_to\_output\_file>]}
\]
to call \gpuourmethod to solve the SDP instance.

\item  If the user has prepared data in sparse SDPA format, created the required config and output files, and computed an appropriate trace bound $\tau$, the user should then execute the command
\[
\texttt{./cuHallar -i <path\_to\_SDPA\_file> --trace\_bound \(\tau\) [-c <path\_to\_config\_file>] [-o <path\_to\_output\_file>]}
\]
to call \gpuourmethod to solve the SDP instance.
\end{itemize}

\paragraph{Runtime reporting.} During each of \gpuourmethod's outer augmented Lagrangian iterations, the current violation in KKT residuals, the current penalty parameter, and the current objective value are all printed for the user to see.

% primal feasibility, the primal dual-gap value, and 

% objective values, residual norms, complementarity, and the current penalty parameter. The option \texttt{output\_path} (or \texttt{-o}) controls where $(Y,p,\theta)$ is saved; $X$ is recovered as $YY^\top$.

\section{Interface}\label{sec:interface}
This section describes the interface and inputs needed to run \gpuourmethod using either HSLR or SDPA format. 

\subsection{Data Input}\label{subsec:data-input}
\ourmethod and \gpuourmethod accept problem data in either HSLR or SDPA format; SDPA inputs are converted in memory to the hybrid layout. HSLR allows each matrix $A_\ell$ to be written as the sum of a sparse matrix and a low rank matrix, i.e.\ (under the convention that $A_0=C$):
\begin{equation}\label{eq:HSLR-rep}
A_\ell = A_\ell^{\mathrm{sp}} \;+\; A_\ell^{\mathrm{lr}}, 
\qquad 
A_\ell^{\mathrm{lr}} = P_\ell D_\ell P_\ell^\top,
\qquad 
\ell=0,1,\dots,m,
\end{equation}
where $A^{\mathrm{sp}} \in\mathbb{S}^n$ is the sparse component of $A_\ell$, and $P_\ell \in \mathbb{R}^{n\times r_\ell}$ and $D_\ell \in\mathbb{S}^{r_{\ell}}$ 
are the factors of the low rank component of $A_\ell$. 
Only the upper triangular part of $A_\ell^{\mathrm{sp}}$ is stored by 
specifying triplets $(i,j,\mathrm{val})$ with $1\le i\le j\le n$ corresponding to its nonzero entries; (duplicate triplets are not allowed). Each column of $P_\ell$ and the corresponding column of $D_\ell$ are provided in a single line of the input file, with the vector for the column of $P_\ell$ separated from the vector for the column of $D_\ell$ by a semicolon `;`. For example, in the maximum stable set formulation of Section~\ref{subsec:max-stable}, the objective matrix $C=-ee^\top$ is dense but has rank one. Representing this matrix in SDPA format would require storing all $\mathcal{O}(n^2)$ nonzero entries, whereas the HSLR format requires only $\mathcal{O}(n)$ storage for the low-rank factor, enabling the solution of much larger instances than would otherwise be feasible.
% The columns of $P_\ell$ and the corresponding columns of $D_\ell$ are stored by columns: each line lists a vector $p\in\mathbb{R}^n$ followed by a scalar $d\in\mathbb{R}$, written as $p_1\,\cdots\,p_n\ ;\ d$, which encodes one column of $P_\ell$ and the corresponding diagonal entry of $D_\ell$; typical use sets $D_\ell$ diagonal. 
The first, second, and third lines of the HSLR file specifies the number of constraints $m$ and matrix size $n$, the right–hand side vector $b\in\mathbb{R}^m$, and the trace bound $\tau>0$, respectively.
Next, each 
$A_\ell$ is specified  by entering  ``\texttt{$\ell$ SP}'' followed by its sparse component description,
and/or by entering ``\texttt{$\ell$ LR}''  followed by its low rank description. Each $A_\ell$ can be entered in any order but, if $A_\ell$ has both sparse and low rank  components, the sparse should precede the low rank one.

The following example is used to illustrate how the input data should be provided.  We first note that all numerical data for matrix entries and vector components are parsed as floating-point values, accepting integer, decimal, and scientific notation allowing for instances like $1e5$. In contrast, the problem dimensions $m$, $n$, and all sparse matrix indices must be provided as integers.

Consider the SDP problem given by
\begin{equation}\label{eq:simex}
\min\ \langle I+ee^\top,X\rangle
\quad\text{s.t.}\quad
\langle 0.5\,I,X\rangle=2,\ \ 
\langle ee^\top,X\rangle=4,\ \ 
%\langle A_3^{\mathrm{sp}}+A_3^{\mathrm{lr}},X\rangle=9,\ \ 
\langle A_3^{\mathrm{sp}}+A_3^{\mathrm{lr}},X\rangle=7,\ \ 
\mathrm{Tr}(X)\le 5,\ \ 
X\succeq 0,
\end{equation}
where $I$ denotes the $4\times4$ identity matrix, $X\in\mathbb{S}^4$, $e:=(1,1,1,1)^\top\in\mathbb{R}^4$, and
% \[
% A_3^{\mathrm{sp}}=\mathrm{diag}(-1,-2,-2,0), 
% \qquad 
% P_3=\begin{bmatrix}1&1\\[0.1em]1&2\\[0.1em]2&1\\[0.1em]1&1\end{bmatrix},
% \qquad 
% D_3=I_2.
% \]
\[
A_3^{\mathrm{sp}}=\begin{bmatrix}
0 & 0 & 0.7 & 0\\
0 & 1 & 0 & -0.5\\
0.7 & 0 & 0 & 0\\
0 & -0.5 & 0 & -1
\end{bmatrix},\qquad
P_3=\begin{bmatrix}
1.0&2.0\\[0.1em]
2.0&1.0\\[0.1em]
1.0&1.0\\[0.1em]
2.0&1.0
\end{bmatrix},\qquad
D_3=\begin{bmatrix}
    1.0 & -0.5\\
    -0.5 & -2.0
\end{bmatrix}.
\]
Since the trace constraint is encoded by the line containing $\tau$, the file records $m=3$ equality constraints and $n=4$, and uses $b=(2,4,7)^\top$. 

The complete HSLR file\footnote{Line spacing is added below for readability but is optional.} is:
\begin{lstlisting}[caption={HSLR file for Problem \eqref{eq:simex}.}]
# m n
3 4
# b vector
2.0 4.0 7.0
# Trace bound
5.0

# Matrix 0: C = I + ee^T
0 SP
1 1 1.0
2 2 1.0
3 3 1.0
4 4 1.0
0 LR
1.0 1.0 1.0 1.0 ; 1.0

# Matrix 1: 0.5 * I
1 SP
1 1 0.5
2 2 0.5
3 3 0.5
4 4 0.5

# Matrix 2: ee^T
2 LR
1.0 1.0 1.0 1.0 ; 1.0

# Matrix 3: A3_sp + A3_lr
3 SP
1 3  0.7
2 2  1.0
2 4 -0.5
4 4 -1.0
3 LR
1.0 2.0 1.0 2.0 ;  1.0   -0.5
2.0 1.0 1.0 1.0 ; -0.5   -2.0
\end{lstlisting}

\subsection{Parameter Input}\label{subsec:params}
Users supply options either on the command line or via a text configuration file. Precedence is strict:
\begin{equation*}
\text{command line} \;\succ\; \text{configuration file} \;\succ\; \text{compiled defaults}.
\end{equation*}
When both are present, the value given on the command line overrides the value in the configuration file, and any value not specified by the user is taken from the solver’s compiled defaults. Paths may be absolute or relative to the current working directory at invocation. Options and their defaults are listed in Table~\ref{tab:options} and grouped as Basic, Intermediate, and Advanced to mirror the solver’s control flow.

The command line accepts short flags for core I/O (\texttt{-i}, \texttt{-p}, \texttt{-d}, \texttt{-c}) and long GNU--style flags for all parameters (e.g., \texttt{--eps\_gap 1e-6}). Boolean flags appear without a value when enabled (e.g., \texttt{--run\_tests}). Numeric quantities use SI units: \texttt{time\_limit} is in seconds; tolerances such as \texttt{eps\_gap} and \texttt{eps\_pfeas} are dimensionless; iteration limits are integers; penalty and Lipschitz parameters are real and must be positive. For SDPA inputs, a trace bound must be provided by the user via \texttt{--trace\_bound} unless encoded in the hybrid HSLR file; for HSLR inputs, the trace bound is read from the file header.

A configuration file is plain text. Each nonempty line assigns a single option by either
\[
\texttt{key = value}\quad\text{or}\quad\texttt{key\ value},
\]
with optional surrounding whitespace. Blank lines are ignored. Lines beginning with \texttt{\#} are treated as comments and ignored. Keys must match the option names in Table~\ref{tab:options}. Values follow the same typing and units as on the command line.

\paragraph{Examples.} The following two invocations are equivalent due to precedence:
\begin{lstlisting}
# Using a configuration file and overriding one option on the command line
$ cuHallar -i model.hslr -c options.cfg --L_inc_fista 3.1

# A minimal options.cfg (values not listed fall back to defaults)
# I/O
input_path = model.hslr
output_path = out.csv
# Stopping criteria
eps_gap = 1e-5
eps_pfeas = 1e-5
time_limit = 3600
# Penalty schedule
beta0 = 10.0
beta_inc = 1.1
beta_min = 10.0
beta_max = 1e11
# ADAP-FISTA
maxiter_fista = 10000
L0_fista = 1.0
L_inc_fista = 2.0
mu_fista = 0.5
chi_fista = 1e-4
sigma_fista = 0.3
err_tol_fista = 1e-8
# AIPP and HLR
maxiter_aipp = 5
lam0_aipp = 0.1
maxiter_hlr = 10
maxiter_hallar = 10000
# Scaling and verbosity
scale_A = 1.0
scale_C = 1.0
verbosity = 1
\end{lstlisting}

In practice, place persistent choices in the configuration file and use the command line to override a small number of run–specific parameters (e.g., \texttt{--trace\_bound} for SDPA input or a one–off change to \texttt{--L\_inc\_fista}). This keeps runs reproducible while preserving exact control through the precedence rule.

\paragraph{Basic Settings.}
These options control the fundamental behavior of the solver. The user must specify the input file path (\texttt{-i}), the primal output file path (\texttt{-p}), the dual output file path (\texttt{-d}), and an optional configuration file (\texttt{-c}). The main termination criteria are also basic settings: the maximum number of outer iterations (\texttt{maxiter\_hallar}), the relative duality gap tolerance (\texttt{eps\_gap}), the primal feasibility tolerance (\texttt{eps\_pfeas}), and the maximum runtime (\texttt{time\_limit}). In the case \texttt{time\_limit} is reached, \gpuourmethod returns the last known iteration information. Finally, the verbosity level (\texttt{verbosity}) controls the amount of information printed to the console. In the case of verbosity equal zero, no output is given to allow for users to implement the models within their own subroutines. 

\paragraph{Intermediate Settings.}
This group of parameters allows for fine-tuning the solver's scaling and penalty updates. The \texttt{scale\_A} and \texttt{scale\_C} options apply uniform scaling to the constraint and cost matrices, respectively, which can improve numerical stability. More details can be found in Subsection~\ref{subsection scaling}. The penalty parameter $\beta$ is controlled by its initial value, \texttt{beta0}, its increment factor \texttt{beta\_inc}, and its lower and upper bounds \texttt{beta\_min} and \texttt{beta\_max}. This category also includes controls for the main inner loop, such as the maximum number of iterations for the Hybrid Low-Rank subroutine which consists of an ADAP-AIPP call + potential Frank-Wolfe steps (\texttt{maxiter\_hlr}) and the initial $\lambda_0$ parameter for the AIPP subroutine (\texttt{lam0\_aipp}).

\paragraph{Advanced Settings.}
These parameters are intended for expert users who wish to control the low-level behavior of the innermost subroutines. This group includes all parameters for the ADAP-FISTA method, such as its iteration limit (\texttt{maxiter\_fista}), step-size parameters (\texttt{mu\_fista}, \texttt{chi\_fista}, \texttt{sigma\_fista}), and Lipschitz constant controls (\texttt{L0\_fista}, \texttt{L\_inc\_fista}). It also contains the tolerances for FISTA (\texttt{err\_tol\_fista}) and the eigenvalue solvers (\texttt{eps\_eig}, \texttt{err\_tol\_eig}), as well as the AIPP iteration limit (\texttt{maxiter\_aipp}). Adjusting these values can impact the trade-off between solution speed and accuracy, but they are generally left at their default values. 

\begin{table}[H]
\centering
\renewcommand{\arraystretch}{1.2}
\begin{tabular}{|l|l|p{8.5cm}|}
\hline
\textbf{Option} & \textbf{Default Value} & \textbf{Description} \\
\hline
\multicolumn{3}{|c|}{\textbf{Input / Output}} \\
\hline
\texttt{-i} & \texttt{none (required)} & Path to the input file in HSLR format; SDPA \texttt{.dat-s} also accepted. \\
\texttt{-p} & \texttt{"primal\_out.txt"} & Path for the output file containing the primal solution. \\
\texttt{-d} & \texttt{"dual\_out.txt"} & Path for the output file containing the dual solution. \\
\texttt{-c} & \texttt{""} & Path to a configuration file to load options. \\
\texttt{--initial\_solution}, \texttt{-w} & \texttt{""} & Path to CSV with dense $Y_0\in\mathbb{R}^{n\times r}$ (no header) for primal warm start; if empty, a default $Y_0$ is used. \\
\texttt{--run\_tests} & \texttt{false (flag)} & Run test routine with example instances. \\
\hline
\multicolumn{3}{|c|}{\textbf{FISTA Parameters}} \\
\hline
\texttt{--maxiter\_fista} & \texttt{1e4} & Maximum number of ADAP-FISTA iterations. \\
\texttt{--mu\_fista} & \texttt{0.5} & FISTA parameter $\mu$. \\
\texttt{--chi\_fista} & \texttt{1e-4} & FISTA parameter $\chi$. \\
\texttt{--L0\_fista} & \texttt{1.0} & Initial Lipschitz constant for ADAP-FISTA. \\
\texttt{--L\_inc\_fista} & \texttt{2.0} & Lipschitz constant increment factor. \\
\texttt{--sigma\_fista} & \texttt{0.3} & FISTA parameter $\sigma$. \\
\texttt{--err\_tol\_fista} & \texttt{1e-8} & Error tolerance for ADAP-FISTA. \\
\hline
\multicolumn{3}{|c|}{\textbf{AIPP Parameters}} \\
\hline
\texttt{--maxiter\_aipp} & \texttt{5} & Maximum number of AIPP iterations. \\
\texttt{--lam0\_aipp} & \texttt{0.1} & AIPP initial parameter $\lambda_0$. \\
\hline
\multicolumn{3}{|c|}{\textbf{Hybrid Low-Rank \& Hallar}} \\
\hline
\texttt{--maxiter\_hlr} & \texttt{10} & Maximum iterations for the hybrid low-rank method. \\
\texttt{--maxiter\_hallar} & \texttt{1e4} & Maximum number of outer \ourmethod iterations. \\
\hline
\multicolumn{3}{|c|}{\textbf{Stopping Criteria}} \\
\hline
\texttt{--eps\_pfeas} & \texttt{1e-5} & Primal feasibility tolerance ($\epsilon_{\text{feas}}$). \\
\texttt{--eps\_gap} & \texttt{1e-5} & Relative duality gap tolerance ($\epsilon_{\text{gap}}$). \\
\hline
\multicolumn{3}{|c|}{\textbf{Penalty Parameters}} \\
\hline
\texttt{--beta0} & \texttt{10.0} & Initial penalty parameter $\beta_0$. \\
\texttt{--beta\_inc} & \texttt{1.1} & Increment factor for $\beta$. \\
\texttt{--beta\_min} & \texttt{10.0} & Minimum value for $\beta$. \\
\texttt{--beta\_max} & \texttt{1e11} & Maximum value for $\beta$. \\
\hline
\multicolumn{3}{|c|}{\textbf{Scaling}} \\
\hline
\texttt{--scale\_A $(\tau_a)$} & \texttt{1.0} & Positive scaling factor for constraint matrices. \\
\texttt{--scale\_C $(\tau_c)$} & \texttt{1.0} & Positive scaling factor for the cost matrix. \\
%\texttt{--trace\_bound $(\tau)$} & \texttt{1.0} & $\tau$ value for the trace constraint. Used when passing a sparse SDPA format file. When using HSLR format, the trace bound is passed inside the input file. \\
\hline
\multicolumn{3}{|c|}{\textbf{Miscellaneous}} \\
\hline
\texttt{--verbosity} & \texttt{1} & Verbosity level (0: silent, 1: summary steps, 2: detailed, 3: debug). \\
\texttt{--time\_limit} & \texttt{3600.0} & Time limit in seconds. \\
\hline
\end{tabular}
\caption{Parameter options grouped by category.}
\label{tab:options}
\end{table}

\subsection{Scaling}\label{subsection scaling}

This subsection explains the roles of the scaling parameters \texttt{scale\_A} and \texttt{scale\_C},
which are denoted by $\tau_a$ and $\tau_c$ in the discussion below.
%\gpuourmethod and \ourmethod support simultaneous scaling of the objective $C$ and the constraint map $\mathcal{A}$.
Recall that the original problem  \eqref{eq:sdp-primal} has the trace constraint $\mathrm{Tr}(X)\le \tau$ and that its optimal value is denoted by $v_*$. Let
\begin{equation}\label{eq:scaled-data}
\tilde C \defi \tau_c C,\qquad 
\mathcal{\tilde A} \defi \tau_a \mathcal A,\qquad 
\tilde b \defi \frac{\tau_a}{\tau} b\,,
\end{equation}
and define the primal and dual scaled variables
\begin{equation}\label{eq:scaled-vars}
\tilde X \defi \frac{1}{\tau} X,\qquad 
\tilde p \defi \frac{\tau_c}{\tau_a}\, p,\qquad 
\tilde \theta \defi \tau_c \theta,\qquad 
\tilde S \defi \tau_c S.
\end{equation}
Then, \eqref{eq:sdp-primal} and \eqref{eq:sdp-dual} are equivalent to the following scaled primal-dual pair  SDPs:
\begin{align}
%\label{eq:scaled-primal}
\tilde v_* &= \min\{\ \tilde C \bullet \tilde X \ :\ \mathcal{\tilde A}(\tilde X)=\tilde b,\ \mathrm{Tr}(\tilde X)\le 1,\ \tilde X\succeq 0\ \} \notag \\
\label{eq:scaled-dual}
&= \max\{\ -\tilde b^{\top}\tilde p - \tilde \theta \ :\ \tilde S=\tilde C+\mathcal{\tilde A}^{*}(\tilde p)+\tilde \theta I \succeq 0,\ \tilde \theta\ge 0\ \},
\end{align}
and there holds
\begin{equation}\label{eq:value-relation}
v_* = \frac{\tau}{\tau_c}\,\tilde v_*.
\end{equation}
The augmented–Lagrangian subproblem associated with \eqref{eq:sdp-primal},
\begin{equation}\label{eq:AL-original}
\min\Big\{\, C\bullet X + \langle p,\,\mathcal A X - b\rangle + \frac{\beta}{2}\,\|\mathcal A X - b\|^2 \ :\ \mathrm{Tr}(X)\le \tau,\ X\succeq 0 \Big\},
\end{equation}
transforms, after multiplying the objective by $\tau_c/ \tau$ and using \eqref{eq:scaled-data}–\eqref{eq:scaled-vars}, into
\begin{equation}\label{eq:AL-scaled}
\min\Big\{\, \tilde C\bullet \tilde X + \langle \tilde p,\,\mathcal{\tilde A}\tilde X - \tilde b\rangle + \frac{\tilde \beta}{2}\,\|\mathcal{\tilde A}\tilde X - \tilde b\|^2 \ :\ \mathrm{Tr}(\tilde X)\le 1,\ \tilde X\succeq 0 \Big\},
\end{equation}
with penalty coupling
\begin{equation}\label{eq:penalty-scaling}
\tilde \beta = \frac{\tau\tau_c}{\tau_a^{2}}\,\beta.
\end{equation}
Relations \eqref{eq:scaled-data}–\eqref{eq:penalty-scaling} make explicit how choices of $(\tau_c,\tau_a)$ relocate magnitude across the objective, constraints, and trace bound.

\section{Output and Interpretation}\label{Output (Diego)}

This section describes the output produced by \gpuourmethod and \ourmethod, including the criteria used for termination, the information reported in the solver logs (console output), and the structure of the solution file.

\subsection{Termination Criteria}
The solvers implement an augmented Lagrangian framework that seeks to satisfy the optimality conditions for the primal-dual pair \eqref{eq:sdp-primal} and \eqref{eq:sdp-dual}. The algorithm terminates when the normalized residuals associated with these conditions fall below specified tolerances.

For user-defined tolerances $\epsilon_{\text{feas}}$ (\texttt{eps\_pfeas}) and $\epsilon_{\text{gap}}$ (\texttt{eps\_gap}), \gpuourmethod stops successfully when an iterate $(X, p, \theta)\in\Delta^n\times\mathbb R^{m}\times\R_{+}$, satisfies the following two conditions:

\begin{enumerate}
    \item \textbf{Primal Feasibility:} The relative error in the equality constraints must be sufficiently small.
    \begin{equation}\label{eq:term_pfeas}
        \frac{\|\mathcal{A}(X)-b\|_2}{1+\|b\|_1} \leq \epsilon_{\text{feas}}.
    \end{equation}

    \item \textbf{Relative Duality Gap:} The difference between the primal objective value ($\mathrm{pval} = C \bullet X$) and the dual objective value ($\mathrm{dval} = -b^{\top}p-\tau\theta$) must be relatively small.
    \begin{equation}\label{eq:term_gap}
        \frac{|\mathrm{pval} - \mathrm{dval}|}{1+|\mathrm{pval}|+|\mathrm{dval}|} \leq \epsilon_{\text{gap}}.
    \end{equation}

    % \item \textbf{Dual Feasibility (Approximate):} The dual slack matrix $S = C+\mathcal A^{*}(p)+\theta I$ must be approximately positive semidefinite.
    % \begin{equation}\label{eq:term_dfeas}
    % \frac{\max\{0,-\lambda_{\min}(S)\}}{1+\|C\|_1} \leq \epsilon_{\text{gap}},
    % \end{equation}
    % where $\lambda_{\min}(S)$ denotes the smallest eigenvalue of $S$.
\end{enumerate}

Our returned solution further satisfies the following properties:

\begin{enumerate}[start=3]
    \item \textbf{Exact Dual Feasibility:} The dual slack matrix is constructed as $S = C+\mathcal A^{*}(p)+\theta I$. Hence, the dual feasibility condition $C+\mathcal A^{*}(p)+\theta I - S = 0$ is satisfied exactly.
    \item \textbf{Exact Primal PSD:} The primal matrix is computed in factorized form $X = Y Y^\top$. Hence, $X$ is always positive semidefinite.
    \item \textbf{Exact Dual PSD:} The choice of dual variable $\theta = \max\{0, -\lambda_{\min}(C + \mathcal A^*(p))\}$ guarantees that $S = C+\mathcal A^{*}(p)+\theta I$ is always positive semidefinite.
\end{enumerate}

If \gpuourmethod reaches the iteration limit (\texttt{maxiter\_hallar}) or the time limit (\texttt{time\_limit}) before satisfying these criteria, it terminates and returns the best solution found so far.

\subsection{Console Output (Solver Log)}
During execution, \gpuourmethod prints a log to the console detailing the progress of the algorithm. The level of detail is controlled by the \texttt{verbosity} parameter. At the default level, the log provides a header summarizing the parameter settings and the problem dimensions, followed by a table reporting the status of each outer Augmented Lagrangian (AL) iteration.

A sample console output is shown below.

\begin{lstlisting}[caption={Sample console output from cuHALLaR.}]
---------- Basic Settings ------------------
input_path = examples/mc_3.dat-s
output_path = out.txt
...
Reading SDPA file: examples/mc_3.dat-s
Problem dimensions:
  - Matrix size: 3000 x 3000
  - Number of constraints: 216172
  - Trace bound: 51601.0

Solving SDP problem with GPU acceleration...

##########################################################################
  #   rank        gap    feas    pval    dval    pnlty   steps
  0    1          -     2.9e-03   9.690e-06    NaN    1.0e+01 A
  1    1        NaN     2.9e-03   8.201e-06   2.500e-03   1.0e+01  A
  2    1        NaN     2.9e-03   6.903e-06   6.250e-03   1.0e+01  A
...
 41    3        9.5e-06   2.1e-08   8.357e-02   8.357e-02   6.0e+03  A
 42    3        8.8e-06   1.3e-08   8.357e-02   8.357e-02   6.6e+03
Final Results
Primal Obj              = 0.08356806847402057
Dual Obj                = 0.08356659006982121
PD Gap                  = 8.844561680506419e-6
Primal infeasibility      = 1.3353696237066644e-8

#ADAP FISTA Calls: 44
#ACG Iterations: 262
#FW Calls: 2
Primal val unscaled = 4312.195901327936
Run time = 2.718115 seconds
Writing output
Output written to primal_out.txt and dual_out.txt.
\end{lstlisting}

The columns in the iteration table are interpreted as follows:
\begin{itemize}
    \item \texttt{\#}: The outer AL iteration count.
    \item \texttt{rank}: The rank $r$ of the current primal iterate $X=YY^\top$, where $Y \in \mathbb{R}^{n \times r}$.
    \item \texttt{gap}: The current relative duality gap, computed according to \eqref{eq:term_gap}.
    \item \texttt{feas}: The current primal feasibility residual, computed according to \eqref{eq:term_pfeas}.
    \item \texttt{pval}: The current primal objective value $C \bullet X$.
    \item \texttt{dval}: The current dual objective value $-b^{\top}p-\tau\theta$.
    \item \texttt{pnlty}: The current value of the penalty parameter $\beta$.
    \item \texttt{steps}: Indicates the type of inner-loop steps taken during the iteration. `A' denotes a call to the ADAP-AIPP subroutine, and `F' denotes a Frank--Wolfe step.
\end{itemize}
The solution \texttt{rank} typically increases after an `F' step, reflecting the addition of a new rank-one component (atom) to the factorization. The penalty parameter \texttt{pnlty} is adjusted adaptively; it generally increases to enforce feasibility but may decrease if subproblem residuals permit.

Upon termination, the log reports the final objective values, gap, and feasibility, followed by statistics on the total number of calls to subroutines (AIPP, FISTA, FW), the unscaled primal objective value (if scaling was applied; see Subsection~\ref{subsection scaling}), and the total runtime.

\subsection{Solution File Output}

\gpuourmethod saves the primal and dual solutions to separate files, whose paths are specified by the user. The final low-rank primal factor $Y \in \mathbb{R}^{n \times r}$ is saved to the path given by \texttt{--primal\_output\_path} (or \texttt{-p}). The file is formatted as a standard comma-separated value (CSV) text file without a header; each of its $r$ columns corresponds to a column of $Y$. The full primal solution matrix is recovered as $X = YY^\top$.

The final dual variables $(p, \theta)$ are saved to a path specified by the flag \texttt{--dual\_output\_path}. The file is formatted as a single comma-separated value (CSV) line; the first field contains the scalar $\theta \ge 0$, and the subsequent $m$ fields contain the components of the vector $p \in \mathbb{R}^m$. The dual slack matrix $S$ is not saved, as it is uniquely determined by the relation $S = C + \mathcal{A}^*(p) + \theta I$ and is guaranteed by construction to be positive semidefinite.

Consider a problem with $n=4$, $m=3$, which terminates with a rank-$2$ solution. The output files would appear as follows.
\begin{lstlisting}[caption={Example primal output file for a rank-2 solution ($Y \in \mathbb{R}^{4\times 2}$).}, label={lst:primal_out}]
# File specified by --primal_output_path out_Y.csv
0.8561,-0.0152
-0.0152,0.9998
-0.5163,0.0021
0.1005,-0.1009
\end{lstlisting}
\begin{lstlisting}[caption={Example dual output file for $m=3$.}, label={lst:dual_out}]
# File specified by --dual_output_path out_p.csv
0.5873,-0.5873,3.4121,-1.2345
\end{lstlisting}
All CSV fields are comma–separated with no embedded whitespace; numeric fields are written in floating–point format. Parsing is unambiguous: $Y$ is read from the $n\times r$ data block in the primal file, and $(\theta,p)$ are read from the single comma-separated line in the dual file.

\appendix

\section{Additional Examples of HSLR Format}\label{sec:examples}

This section displays how to construct HSLR format for several structured SDPs such as the SDP relaxations of the Matrix Completion and Maximum Stable Set problems. The main quantities that are needed for HSLR format are the dimension pair $(m,n)$, the cost matrix $C$, the data matrices $A_{\ell}$, $\ell=1,\ldots m$, the right-hand side vector $b$, and the tracebound $\tau$. By convention we consider $C$ to be matrix $0$, i.e., $C=A_0$. Recall from Subsection~\ref{subsec:data-input} that the matrices $A_{\ell}$ are assumed to have the structure
\begin{equation}\label{eq:HSLR-rep-Stable}
A_\ell = A_\ell^{\mathrm{sp}} + A_\ell^{\mathrm{lr}}, 
\qquad 
A_\ell^{\mathrm{lr}} = P_\ell D_\ell P_\ell^\top,
\qquad 
\ell=0,1,\dots,m,
\end{equation}
where $A^{\mathrm{sp}} \in\mathbb{S}^n$ is the sparse component of $A_\ell$, and $P_\ell \in \mathbb{R}^{n\times r_\ell}$ and $D_\ell \in\mathbb{S}^{r_{\ell}}$ 
are the factors of the low rank component of $A_\ell$. Subsections~\ref{subsec:matrix-completion} and \ref{subsec:max-stable} below display how the user should prepare HSLR format for the Matrix Completion and Maximum Stable Set SDP relaxations, respectively.

\subsection{Matrix Completion}\label{subsec:matrix-completion}

Given integers $n_2 \ge n_1 \ge 1$, the goal of the matrix completion problem is to recover a low-rank matrix $M \in \mathbb{R}^{n_1 \times n_2}$ by observing a subset of its entries $\{M_{ij} : (i,j) \in \Omega\}$. A standard convex relaxation replaces the rank function with the nuclear norm:
\[
\min_{Y \in \mathbb{R}^{n_1 \times n_2}} \ \{ \ \|Y\|_* \ :\ Y_{ij} = M_{ij}, \ \forall (i,j) \in \Omega \ \}.
\]
Using the semidefinite representation of the nuclear norm, this optimization problem is equivalent to the following SDP
\begin{align}
\label{eq:matcomp-sdp}
\min_{X \in \mathbb{S}^{n_1+n_2}}
\left\{ \tfrac12 \, \mathrm{Tr}(X) \ : \
X = \begin{pmatrix} W_1 & Y \\[0.2em] Y^\top & W_2 \end{pmatrix} \succeq 0, \ 
Y_{ij} = M_{ij} \ \forall (i,j) \in \Omega \right\}.
\end{align}
\gpuourmethod solves the formulation in \eqref{eq:matcomp-sdp}.

In the above formulation, the size of the matrix variable $X$ is $n = n_1 + n_2$ and the number of constraints is $m=|\Omega|$. The right hand side vector $b$ is a vector in $\mathbb R^{m}$ and the $\ell$-th component of $b$ is just $M_{ij}$ where $(i,j)$ is the $\ell$-th index pair in $\Omega$. To compute a suitable tracebound $\tau$, we generate $\hat Y \in \mathbb R^{n_1\times n_2}$ so that $\hat Y_{ij}=M_{ij}$ for indices $(i,j)\in \Omega$ and $\hat Y_{ij}=0$ for all other indices. The tracebound $\tau$ is then computed to be $2 \sqrt{n_1} \|\hat Y\|_F$.

% which is the trivial completion, i.e.,

%$M_{i$b_{ell}$

All data matrices which encode the SDP in \eqref{eq:matcomp-sdp} are sparse. Clearly, $C=A_0=0.5I$. All constraint matrices $A_{\ell}$, $l=1,\ldots m$, have exactly 2 nonzero entries. To see this, consider the constraint $Y_{ij}=M_{ij}$ where $1\leq i\leq n_1$ and $1\leq j\leq n_2$. The $A_{\ell}$ matrix corresponding to this constraint then simply takes value $0.5$ in its positions $(i,n_1+j)$ and $(n_1+j, i)$
and value 0 elsewhere
since this enforces that $A_{\ell}\bullet X=Y_{ij}$. The following example illustrates how to prepare a HSLR data file for a small Matrix Completion problem where $n_1=2$ and $n_2=2$.

% We follow the convention in \cite{monteiro2024low} and set $\tau=3n$.

% 
% If $\Omega = \{(i,j)\}$ with $1 \le i \le n_1$ and $1 \le j \le n_2$, the associated constraint matrix $A_\ell$ enforces $Y_{ij} = M_{ij}$. In $X$, this entry is at $(i, n_1 + j)$, so the matrix $A_\ell$ has the value $0.5$ in its positions $(i, n_1+j)$ and $(n_1+j, i)$, ensuring $A_\ell \bullet X = Y_{ij}$. The right-hand side is $b_\ell = M_{ij}$. 

% The number of constraints is $m = |\Omega|$.

% In HSLR form, all $A_\ell$ matrices are purely sparse:
% \[
% \text{Objective: } C = \tfrac12 I, \quad A_0^{\mathrm{sp}} = \tfrac12 I,\quad A_0^{\mathrm{lr}} = 0.
% \]

\medskip
\noindent\emph{Example:} Suppose that $n_1 = 2$, $n_2 = 2$, $\Omega = \{(1,1), (2,2)\}$, and $M_{11} = 5$ and $M_{22} = 3$. Then $m=2$ and $n=4$.

% Let 
% \[\hat Y=\]

% $\hat Y_{11}=5$

The cost matrix is $C=A_0=0.5I$. It is easy to see that the constraint matrices $A_1$ and $A_2$ and the right-side vector $b$ are:
\begin{align*}
A_1=\begin{bmatrix}
0 & 0 & 1 & 0\\
0 & 0 & 0 & 0\\
1 & 0 & 0 & 0\\
0 & 0 & 0 & 0\\
\end{bmatrix}, \quad A_2=\begin{bmatrix}
0 & 0 & 0 & 0\\
0 & 0 & 0 & 1\\
0& 0 & 0 & 0\\
0 & 1 & 0 & 0\\
\end{bmatrix}, \quad b=\begin{bmatrix}
5\\
3\\
\end{bmatrix}.
\end{align*}
To compute the trace bound, let
$\hat Y=\begin{bmatrix}
5 & 0\\
0 & 3
\end{bmatrix}$. The trace bound $\tau$ is then computed to be $\tau=2\sqrt{2}\|\hat Y\|_{F}\approx 16.50$. The HSLR data file corresponding to this example is:

\begin{lstlisting}[caption={HSLR data file for a small Matrix Completion problem.}]
# m n
2 4
# b vector
5.0 3.0
# Trace bound
16.50

# Matrix 0: C = 0.5*I
0 SP
1 1 0.5
2 2 0.5    
3 3 0.5
4 4 0.5

# Matrix 1: A1_sp
1 SP
1 3 0.5

# Matrix 2: A2_sp 
2 SP
2 4 0.5
\end{lstlisting}

Several remarks about the above HSLR data file are now given. The entries ``2 4'' refer to the dimension pair $(m,n)$ while the entries ``5.0 3.0'' refer to $b_1$ and $b_2$, respectively. The trace bound $\tau$ is 16.50. The line ``Matrix 0: C = 0.5*I'' refers to the fact that the cost matrix is just $0.5*I$. The line ``0 SP'' just means that this matrix is the $0$-th matrix and it is a sparse matrix. The entries ``1 1 0.5'', $\ldots$ ``4 4 0.5'' mean that $C_{11}=0.5, \ldots C_{44}=0.5$. The line ``Matrix 1: A\_sp'' means that we are now writing the first constraint matrix $A_1$. The line ``1 SP'' means that this is the first constraint matrix and it is a sparse matrix. The line ``1 3 0.5'' means that $(A_{1})_{13}=0.5$. Note that this is the only entry that needs to be encoded for $A_1$ since it is the only nonzero entry of the upper triangular part of $A_1$. Finally, the line ``Matrix 2: A2\_sp'' means that we are now writing the second constraint matrix $A_2$. The line ``2 SP'' means that this is the second constraint matrix and it is a sparse matrix. The line ``2 4 0.5'' means that $(A_{2})_{24}=0.5$. 
% \begin{lstlisting}[caption={HSLR data file for a small Matrix Completion problem.}]
% 2 4
% # nnz(A_l^sp) for l=0,...,m
% 4 2 2
% # rank(A_l^lr) for l=0,...,m
% 0 0 0
% # b vector
% 5 3
% # Trace bound
% 12

% # Matrix 0: C = 0.5 * I
% 0 SP
% 1 1 0.5
% 2 2 0.5
% 3 3 0.5
% 4 4 0.5

% # Matrix 1: X_13 = 5
% 1 SP
% 1 3 0.5
% 3 1 0.5

% # Matrix 2: X_24 = 3
% 2 SP
% 2 4 0.5
% 4 2 0.5
% \end{lstlisting}

\subsection{Maximum Stable Set}\label{subsec:max-stable}

% \jacob{Below is a revised subsection, please feel free to make any changes Arnesh.}

For a given undirected graph $G=(V,E)$ with $|V|=n$ vertices and $|E|$ edges, the stability number $\alpha(G)$ is the maximum size of a stable set (a subset of vertices where no two vertices are adjacent). The Lov\'asz $\vartheta$-function provides an upper bound on $\alpha(G)$ and is defined via the following semidefinite program:
\begin{align}\label{eq:lovasz}
\vartheta(G) = \max_{X \in \mathbb{S}^n}\ \{\ J \bullet X \ :\ \mathrm{Tr}(X)\leq1,\ X_{ij}=0 \ \forall \{i,j\}\in E,\ X\succeq 0\ \},
\end{align}
where $J=ee^\top$ is the $n \times n$ matrix of all ones, and $e \in \mathbb{R}^n$ is the vector of all ones.

In this formulation, the size of the matrix variable $X$ is $n=|V|$ and the number of constraints is $m = |E|$. 
% The trace constraint $\mathrm{Tr}(X)=1$ constraint is not included in $m$ since we treat this constraint as our implicit trace bound constraint $\mathrm{Tr}(X)\leq \tau$. 
The trace bound $\tau$ is thus set to be 1. The right hand side vector $b\in \mathbb R^{m}$ is simply the vector of all zeros.

To adapt this formulation to the minimization format required by our HSLR format, we minimize the negative of the objective function, setting the cost matrix to $C = -J$. 
% We now detail the structure of the data matrices for the HSLR format. 
The cost matrix $C=A_0=-J$ is dense but has rank one. It is represented efficiently using only its low-rank component:
\[
A_0^{\mathrm{sp}} = 0, \quad A_0^{\mathrm{lr}} = P_0 D_0 P_0^\top, \quad \text{where } P_0 = e \in \mathbb{R}^{n\times 1} \text{ and } D_0 = [-1] \in \mathbb{S}^1.
\]
In HSLR format, only the components $D_0$ and $P_0$ need to be specified. This representation is significantly more storage-efficient than storing the fully dense matrix $-J$. 

All constraint matrices $A_{\ell}$, $\ell=1,\ldots m$, are sparse and have exactly two nonzero entries. To see this, consider the constraint $X_{ij}=0$ where $\{i,j\}$ is an edge. The $A_{\ell}$ matrix corresponding to this constraint is then just $A_\ell = 0.5(E_{ij} + E_{ji})$, where $E_{ij}$ denotes the matrix with $1$ in position $(i,j)$ and zeros elsewhere. This construction enforces that $A_{\ell}\bullet X=X_{ij}$. Since HSLR format only requires the user to specify the nonzero components of the upper triangular part of a sparse matrix, the user only needs to specify one entry for each of the constraint matrices $A_{\ell}$, $\ell=1,\ldots m$.  The following example illustrates how to prepare a HSLR data file for a small Maximum Stable Set problem where $n=4$ and $m=4$.

% For an edge $\{i,j\}$, the constraint $X_{ij}=0$ is enforced by setting $A_\ell = 0.5(E_{ij} + E_{ji})$, where $E_{ij}$ denotes the matrix with $1$ in position $(i,j)$ and zeros elsewhere. This ensures that $A_\ell \bullet X = X_{ij}$.

% This representation only requires $\mathcal{O}(n)$ storage compared to $\mathcal{O}(n^2)$ for the dense matrix.

% The constraints $X_{ij}=0$ can all be represented

% include the trace constraint $\mathrm{Tr}(X)=1$ and the edge constraints $X_{ij}=0$.

% and hence set 

% The trace bound is set to $\tau=1$. S

% ince $\mathrm{Tr}(X)=1$ is explicitly included in the constraints $\mathcal{A}(X)=b$, the condition $\mathrm{Tr}(X)\le \tau$ in \eqref{eq:sdp-primal} is satisfied.

% The constraints are all sparse. The first constraint (indexed $\ell=1$) corresponds to $\mathrm{Tr}(X)=1$. The constraint matrix is $A_1 = I$, and the right-hand side is $b_1=1$.

% The subsequent constraints, indexed from $\ell=2$ to $m$, correspond to the edges $\{i,j\} \in E$. For an edge $\{i,j\}$, the constraint $X_{ij}=0$ is enforced by setting $A_\ell = 0.5(E_{ij} + E_{ji})$, where $E_{ij}$ denotes the matrix with $1$ in position $(i,j)$ and zeros elsewhere. This ensures $A_\ell \bullet X = X_{ij}$. The matrix $A_\ell$ is sparse, having only two nonzero entries, and the right-hand side is $b_\ell=0$.

\medskip
\noindent\emph{Example:} Consider the 4-cycle graph $C_4$, with $V=\{1,2,3,4\}$ and $E=\{\{1,2\}, \{2,3\}, \{3,4\}, \{1,4\}\}$. We have $|V|=4$ and $|E|=4$, so $n=4$ and $m=4$. The trace bound is $\tau=1.0$. The cost matrix $C=A_0=-J$. Its low-rank factorization uses:
\[
P_0 = (1, 1, 1, 1)^\top, \quad D_0 = [-1].
\]
The constraint matrices $A_1, \ldots, A_4$ and the right-hand side vector $b$ are:
\begin{equation}
   A_1=\begin{bmatrix}
0 & 0.5 & 0 & 0\\
0.5 & 0 & 0 & 0\\
0 & 0 & 0 & 0\\
0 & 0 & 0 & 0\\
\end{bmatrix}, \quad A_2=\begin{bmatrix}
0 & 0 & 0 & 0\\
0 & 0 & 0.5 & 0\\
0 & 0.5 & 0 & 0\\
0 & 0 & 0 & 0\\
\end{bmatrix},
\end{equation}

\begin{equation}
A_3=\begin{bmatrix}
0 & 0 & 0 & 0\\
0 & 0 & 0 & 0\\
0 & 0 & 0 & 0.5\\
0 & 0 & 0.5 & 0\\
\end{bmatrix}, \quad A_4=\begin{bmatrix}
0 & 0 & 0 & 0.5\\
0 & 0 & 0 & 0\\
0 & 0 & 0 & 0\\
0.5 & 0 & 0 & 0\\
\end{bmatrix}, \quad b=\begin{bmatrix}
0\\
0\\
0\\
0
\end{bmatrix}.  
\end{equation}

% \begin{align*}
% A_1=\begin{bmatrix}
% 0 & 0.5 & 0 & 0\\
% 0.5 & 0 & 0 & 0\\
% 0 & 0 & 0 & 0\\
% 0 & 0 & 0 & 0\\
% \end{bmatrix}, \quad A_2=\begin{bmatrix}
% 0 & 0 & 0 & 0\\
% 0 & 0 & 0.5 & 0\\
% 0 & 0.5 & 0 & 0\\
% 0 & 0 & 0 & 0\\
% \end{bmatrix}, \\ A_3=\begin{bmatrix}
% 0 & 0 & 0 & 0\\
% 0 & 0 & 0 & 0\\
% 0 & 0 & 0 & 0.5\\
% 0 & 0 & 0.5 & 0\\
% \end{bmatrix}, \quad A_4=\begin{bmatrix}
% 0 & 0 & 0 & 0.5\\
% 0 & 0 & 0 & 0\\
% 0 & 0 & 0 & 0\\
% 0.5 & 0 & 0 & 0\\
% \end{bmatrix}, \quad b=\begin{bmatrix}
% 0\\
% 0\\
% 0\\
% 0
% \end{bmatrix}.
% \end{align*}

The HSLR data file corresponding to this example is:

\begin{lstlisting}[caption={HSLR data file for the Maximum Stable Set problem on $C_4$.}]
# m n
4 4
# b vector
0.0 0.0 0.0 0.0
# Trace bound
1.0

# Matrix 0: C = -J (Low Rank)
0 LR
1.0 1.0 1.0 1.0 ; -1.0

# Matrix 1: Edge (1,2)
1 SP
1 2 0.5

# Matrix 2: Edge (2,3)
2 SP
2 3 0.5

# Matrix 3: Edge (3,4)
3 SP
3 4 0.5

# Matrix 4: Edge (1,4)
4 SP
1 4 0.5
\end{lstlisting}

Several remarks about the above HSLR data file are now given. The first line ``4 4'' specifies the dimension pair $(m,n)$. The second line specifies the vector $b$ which in this case is just $b=(0,0,0,0)^\top$. The third line sets the trace bound $\tau=1.0$. The block starting with ``0 LR'' defines the cost matrix $A_0=C$. The line ``1.0 1.0 1.0 1.0 ; -1.0'' defines the single column of $P_0$ (the vector $e$) and the corresponding entry in $D_0$ (the scalar $-1$), separated by ``;''. The block starting with ``1 SP'' defines $A_1 $. The line ``1 2 0.5'' specifies that $(A_1)_{12} = 0.5$. The block starting with ``2 SP'' defines $A_2$. The line ``2 3 0.5'' specifies that $(A_2)_{23} = 0.5$. The blocks starting with ``3 SP'' and ``4 SP'' define matrices $A_3$ and $A_4$, respectively, in a similar-like fashion. 

% as described for $A_1$ and $A_2$.

% Since only the upper triangular part is required, this single entry fully defines the sparse matrix $A_2$. Note that for $A_5$ (edge $\{4,1\}$), we input the upper triangular indices ``1 4 0.5''.

\medskip

\section*{Acknowledgments}
This research was supported in part through research cyberinfrastructure resources and services provided by the Partnership for an Advanced Computing Environment (PACE) at Georgia Tech, Atlanta, Georgia, USA.

% \newpage
% \bibliographystyle{siam}
\bibliographystyle{plain}
\bibliography{references}
\end{document}